\documentclass[11pt]{article}
\usepackage{amssymb,amsmath}
\usepackage{graphicx}
\usepackage[latin1]{inputenc}
\newtheorem{dfn}{Definition}[section]
\newtheorem{rmk}{Remark}[section]
\newtheorem{thm}{Theorem}[section]
\newtheorem{thm*}{Theorem}
\newtheorem{cor}{Corollary}[section]
\newtheorem{prop}{Proposition}[section]

\newtheorem{ex}{Example}[section]
\newcommand{\Pf}{{\em Proof}. }
\newcommand{\EPf}{\begin{flushright} $\Box$ \end{flushright}}

\newcommand{\ca}{{\rm ca}}

\begin{document}
\title{Limit of Gaussian and normal curvatures  of surfaces in Riemannian approximation scheme for sub-Riemannian three dimensional manifolds and Gauss-Bonnet theorem}
\author{ 
       Jos\'e M. M. Veloso
         \\ Instituto de Ci\^encias Exatas e Naturais
         \\ Universidade Federal do Par\'a
                  \\ veloso@ufpa.br
                  \\https://orcid.org/0000-0002-8969-7320 }

\maketitle
\begin{abstract}
    The authors Balogh-Tyson-Vecchi  in \cite{BTV} utilize the Riemannian approximations scheme $(\mathbb H^1,<,>_L)$, in the Heisenberg group, introduced by Gromov \cite{G}, to calculate the limits of Gaussian and normal curvatures defined on surfaces of $\mathbb H^1$ when $L\rightarrow\infty$. They show that these limits exist (unlike the limit of Riemannian surface  area form or length form), and they obtain Gauss-Bonnet theorem in $\mathbb H^1$ as limit of Gauss-Bonnet theorems in $(\mathbb H^1,<,>_L)$ when $L$ goes to infinity. This construction was extended by Wang-Wei in \cite{WW} to the affine group and the group of rigid motions of the Minkowski plane. %Another approach was given by  Diniz-Veloso who obtained in \cite{DV} Gaussian curvature for non-horizontal surfaces in sub-Riemannian Heisenberg space $\mathbb H^1$ by a different method, and   a Gauss-Bonnet Theorem was proved. The definition was analogous to Gauss curvature of surfaces in $\mathbb R^3$ with particular normal to surface and Hausdorff measure of area. The image of Gauss map was in the cylinder of radius one. That Gaussian curvatures of \cite{BTV} and \cite{DV} does not coincide was proved in \cite{V}.
    We generalize constructions of \cite{BTV} and \cite{WW} to surfaces in sub-Riemannian three dimensional manifolds following the approach of \cite{V}, and prove analogous Gauss-Bonnet theorem.

%The analogous  rotation surfaces of constant curvature defined in  \cite{DV}  are in another paper \cite{DSV2}.
%In \cite{DV} we introduced a notion of  Gaussian curvature, and here we classify rotation surfaces which are of constant Gaussian curvature. We can study these surfaces rotating a horizontal curve $\gamma$ around the $z$-axis. They are very similar to rotation surfaces of constant curvature in Euclidean space $\mathbb R^3$.

 %The rotation surfaces of constant mean curvature in $\mathbb H^1$ are well known. 
%The mean curvature of a rotation surface $S$ in $\mathbb H^1$ is the curvature of the  curve in $\mathbb R^2$ which is the projection of $\gamma$. The  approach we used to classify the rotation surfaces of constant Gaussian curvature is more simple and geometric, and we use this method also to give the rotation surfaces of mean curvature.

\end{abstract}


\section{Introduction}% In \cite{DV}  Gaussian curvature for non-horizontal surfaces in sub-Riemannian Heisenberg space $\mathbb H^1$ was defined and  a Gauss-Bonnet Theorem was proved. The definition was analogous to Gauss curvature of surfaces in $\mathbb R^3$ with particular normal to surface and Hausdorff measure of area. The image of Gauss map was in the cylinder of radius one.

In  \cite{BTV}, a version of Gauss-Bonnet theorem is given through a limit of Riemannian approximations. This work was extended in \cite{WW} to   Gauss-Bonnet theorems in the affine group and in the group of rigid motions of the Minkowski plane. The but of this paper is to generalize these constructions of Gauss-Bonnet theorem to general three dimensional sub-Riemannian manifolds, following the same approach of \cite{V}.  Another way to get a different Gauss-Bonnet theorem in Heisenberg space is in \cite{DV}.
%Furthermore, we show that the limit $K$ and $k_n$ of Riemannian Gaussian curvatures $K^L$ and normal curvatures $k_n^L$ depend only on one of the two functions that define the geometry of the surface. Also we obtain that $K^\infty d\sigma=d(k_nds)$ and  it is possible to get Gauss-Bonnet theorem of \cite{BTV} applying Stokes theorem without taking limit.

A sub-Riemannian manifold $(M,D,<,>)$ is a manifold $M$, with a two-dimensional  contact distribution $D\subset TM$, and a scalar product $<,>$ on $D$.
%The space $\mathbb H^1$ is a Lie group. We define a distribution $D$ generated by the left invariant vector fields $e_1=\frac{\partial}{\partial x}-\frac 1 2y\frac{\partial}{\partial z}$ and $e_2=\frac{\partial}{\partial y}+\frac 1 2x\frac{\partial}{\partial z}$ on $\mathbb H^1$. 
If $e_1,e_2$ is an orthonormal basis of $D$ consider a one form $\omega$ such that $\omega(D)=0$ and $d\omega(e_1,e_2)=-1$. There is a vector field $e_3$ such that $\omega(e_3)=1$ and $i(e_3)d\omega=0$. Then $e_1,e_2,e_3$ is a basis of $TM$.
%are orthonormal. We complete $e_1,e_2$ to a basis of left invariant vector fields in $\mathbb H^1$ by introducing $e_3=[e_1,e_2]=\frac{\partial}{\partial z}$. %Therefore, if $e^1, e^2, e^3$ are dual forms to $e_1,e_2,e_3$, then the volume element invariant by the group action is $dV=e^1\wedge e^2\wedge e^3$.
In the same way of \cite{BTV}, \cite{CDPT}, \cite{V}, \cite{WW}  we consider the family of $<,>_L$ metrics on $TM$ such that $e_1$, $e_2$, $e_3^L=e_3/\sqrt{L}$ is an orthonormal basis of $(TM,<,>_L)$, and the Levi-Civitta connections $\overline\nabla^L$ on $(TM,<,>_L)$.
%The dual basis is $e^1$, $e^2$, $e^3_L=\sqrt{L}e^3$. Then the Carnot-Caratheodory metric space $\mathbb H^1$ is the limit in the sense of Gromov-Hausdorff of Riemannian metric spaces $(R^3,d_L)$, when $L\rightarrow \infty$.  

%We  define a distance in $\mathbb H^1$ by considering the distance between two points as the infimum of length of curves tangent to $D$ that connect the two points. With this distance $\mathbb H^1$ is a metric space with Hausdorff dimension 4, and the differentiable surfaces have Hausdorff dimension 3.
At points of a surface $S$, where the distribution $D$ does not coincide with $TS$, the intersection $D\cap TS$ has dimension 1, and we obtain a direction that we call \emph{characteristic} at this point of $S$. We suppose that all points of  surface $S$ have this property. 
The vector field  \emph{normal horizontal} $f_1$ is a unitary vector field in $D$ orthogonal to the characteristic direction which we suppose to be globally defined. We define $f^1$ by $f^1(f_1)=1$ and $f^1(TS)=0$. We denote by $f_2$ a unitary vector field in $D\cap TS$ and complete a basis of $TS$ taking $f_3=e_3+Af_1$. If $\alpha$ is the angle between $e_1$ and $f_1$, then $f_1=\cos\alpha\, e_1+\sin\alpha\, e_2$ and $f_2=-\sin\alpha\, e_1+\cos\alpha\, e_2$.
%Given a compact set  $K\subset S$, the 3-dimensional (spherical) Hausdorff measure of $K$ is given by $\int_Ki(f_1)dV$ and the curves in $\mathbb H^1$ transverses to $D$ have (spherical) Hausdorff dimension 2 and its Hausdorff measure is $\int_{\gamma}e^3$. 
An orthonormal basis of $TS$ in $(TM,<,>_L)$ is given by $X_2^L=f_2$, $X_3^L=\frac{1}{\sqrt{L+A^2}}f_3$. The normal vector in the scalar product  $<,>_L$ to $S$ is 
$$X_1^L=\frac{\sqrt{L}}{\sqrt{L+A^2}}f_1-\frac{A}{\sqrt{L+A^2}}e_3^L.$$ 
The covariant derivative $\overline \nabla^L$ projects on $S$ through the normal $X_1^L$ and define a covariant derivative $\nabla^L$ on $S$. The Gaussian curvature of $S$ applying the definition given by $\nabla^L$ is $K_L$ and 
$$K=\lim_{L\rightarrow \infty}K_L=-dA(f_2)-A^2.$$
Similarly, we can define the normal curvature $k_n=\lim_{L\rightarrow\infty}k_n^L$ of a transverse curve $\gamma$ in $S$, where  $k_n^L$ is the normal curvature of $\gamma$  applying  $\nabla^L$. Suppose $\gamma(t)$ is a curve in $S$ such that $\gamma'(t)=x(t)f_2(\gamma(t))+y(t)f_3(\gamma(t))$, where $y(t)\neq 0$ for every $t$. Then 
$$
k_n=\lim_{L\rightarrow\infty} k_n^L=\frac{y}{|y|}A.
$$
Then $K$ and $k_n$ depend only on $A$.
If $\sigma _L$ is the area element os $S$ in the metric defined by scalar product $<,>_L$, then the  Hausdorff measure $d\sigma$ on $S$ as a metric space with the sub-Riemannian distance defined by $(M,D,<,>)$ is given by (cf. \cite{BTV}, \cite{DSV})
$$
d\sigma=\lim_{L\rightarrow \infty }\frac{1}{\sqrt{L}}d\sigma_L=f^2\wedge f^3.
$$
In the same way, if $ds_L$ is the length element on $\gamma$ in the metric defined by scalar product $<,>_L$, then the Hausdorff measure $ds$ for transversal curves in $S$ as a metric space with the sub-Riemannian distance defined by $(M,D,<,>)$ is given by (cf. \cite{BTV}, \cite{DSV})
$$
ds=\lim _{L\rightarrow\infty}\frac{1}{\sqrt{L}}ds_L=\frac{y}{|y|}f^3.
$$

%Without necessity to

To get formulas type Gauss-Bonnet obtained in \cite{BTV}, \cite{WW} it is not necessary to apply limits to Gauss-Bonnet formulas in $(M,<,>_L)$. It follows from Stokes'  theorem:
$$
\int_cK d\sigma=%\int_c d\Omega^3_2(f_2,f_3) f^2\wedge f^3=
-\int_c (dA(f_2)+A^2)f^2\wedge f^3=-\int_c d(Af^3)=-\int_{\partial c}Af^3=-\int_{\partial c}k_nds.
$$

\section{Sub-Riemannian manifolds of dimension three}
Suppose $M$ is a three dimensional manifold, $D\subset TM$ a non-integrable distribution with a scalar product $<,>$ in $D$. Choose   $e_1$,  $e_2$, $e_3$ an orthonormal basis of $D$ (may be on an open set $U$ of $M$) as above. 
%Let be $\omega$ a section of $TM^*$ such that $\omega|_D=0$ and $d\omega(e_1,e_2)=-1$. We can choose $e_3$ such that $\omega(e_3)=1$  and  $i(e_3)d\omega=0$. Then $e_1,e_2,e_3$ is a basis of $TM$. 
We have
$$
[e_1,e_2]=\Sigma_{j=1}^2a_{12}^je_j+e_3,
$$
$$[e_i,e_3]=\Sigma_{j=1}^2a_{i3}^j e_j,
$$
for $i=1,2$.
The dual basis $e^1$, $e^2$, $e^3$ satisfies
$$
\begin{array}{rcl}
de^1&=&-a_{12}^1e^1\wedge e^2-a_{13}^1e^1\wedge e^3-a_{23}^1e^2\wedge e^3\\
de^2&=&-a_{12}^2e^1\wedge e^2-a_{13}^2e^1\wedge e^3-a_{23}^2e^2\wedge e^3\\
de^3&=&-e^1\wedge e^2.
\end{array}$$
Differentiating the last equation we get $a_{13}^1+a_{23}^2=0$.

\section{The approximation scheme by scalar products $<,>_L$}                                                                                                                                                                                                                                                                                                                                     
Consider the scalar product $<,>_L$ in $TM$ such that $e_1^L=e_1$, $e^L_2=e_2$, $e_3^L=e_3/\sqrt{L}$ is an orthonormal basis. The dual basis is $e^1_L=e^1$, $e^2_L=e^2$, $e^3_L=\sqrt{L}e^3$. Then the Carnot-Caratheodory metric space $(M,d)$ is the limit in the sense of Gromov-Hausdorff of Riemannian metric spaces $(M, d_L)$, when $L\rightarrow \infty$.  We consider the Levi-Civitta connection $\overline\nabla^L$ in $(M,<,>_L)$. 
This connection is defined by Koszul formula
\begin{equation}\label{koszul}
2<\overline\nabla^L_{e_i^L}e_j^L,e_k^L>_L=<[e_i^L,e_j^L],e_k^L>_L-<[e_j^L,e_k^L],e_i^L>_L+<[e_k^L,e_i^L],e_j^L>_L
\end{equation}
and if $\overline\nabla^Le_i^L=\Sigma_{i=1}^3\omega_i^{Lj}e_j^L$ then $\omega_i^{Lj}$ satisfies $\omega_i^{Lj}=-\omega_j^{Li}$, $1\leq i,j,k\leq 3$.
\begin{prop}
We have
\begin{equation}\label{omegaL}
\begin{array}{rcl}
\omega_1^{L2}&=&-a_{12}^1e^1_L-a_{12}^2e^2_L+\frac{1}{2\sqrt{L}}(a^1_{23}-a^2_{13}-L)e^3_L;\vspace{0.2cm}\\
\omega_1^{L3}&=&-\frac{1}{\sqrt{L}}a_{13}^1e^1_L-\frac{1}{2\sqrt{L}}(a_{13}^2+a_{23}^1+L)e^2_L;\vspace{0.2cm}\\
\omega_2^{L3}&=&\frac{1}{2\sqrt{L}}(-a_{13}^2-a_{23}^1+L)e^1_L-\frac{1}{\sqrt{L}}a_{23}^2e^2_L.
\end{array}
\end{equation}
\end{prop}
\Pf It is a direct application of formula \ref{koszul}.\EPf

\begin{dfn}\label{limomega}Define 
$$
\begin{array}{rcl}
\omega_1^{2}&=&\lim_{L\rightarrow \infty}\frac{1}{L}\omega_1^{L2}=-\frac{1}{2}e^3;\vspace{0.2cm}\\
\omega_1^{3}&=&\lim_{L\rightarrow \infty}\frac{1}{\sqrt{L}}\omega_1^{L3}=-\frac{1}{2}e^2;\vspace{0.2cm}\\
\omega_2^{3}&=&\lim_{L\rightarrow \infty}\frac{1}{\sqrt{L}}\omega_2^{L3}=\frac{1}{2}e^1.
\end{array}
$$
\end{dfn}
%\begin{rmk}
%With the forms $\omega_i^j$ it is possible to define a covariant derivative $\overline\nabla$ by
%$$
%\overline\nabla e_i=\Sigma_{j=1}^3\omega_i^je_j.
%$$
%\begin{prop}The torsion $\overline T$ of $\overline \nabla$ is given by
%This covariant derivative in general has torsion no null, given by
%$$
%T=-e^1\wedge e^2(a^1_{12}e_1+a^2_{12}e_2)-e^1\wedge e^3%\otimes
%(a_{13}^1e_1+a_{13}^2e_2)-e^2\wedge e^3%\otimes
%(a_{23}^1e_1+a_{23}^2e_2).
%$$

%\end{prop}
%\Pf In fact, applying the definition of torsion we get
%$$
%T(e_1,e_2)=(\omega_2^1e_1+\omega_2^3e_3)(e_1)-(\omega_1^2e_2+\omega_1^3e_3)(e_2)-e_3=\frac 1 2 e_3+\frac 1 2 e_3-e_3=0;
%$$
%$$
%\begin{array}{rcl}
%T(e_1,e_3)&=&(\omega_3^1e_1+\omega_3^2e_2)(e_1)-(\omega_1^2e_2+\omega_1^3e_3)(e_3)-(a_1^1e_1+a_1^2e_2+a_1^3e_3)\\
%&=&-\frac 1 2 e_2+\frac 1 2 e_2+a_1^3e_3-(a_1^1e_1+a_1^2e_2+a_1^3e_3)\\
%&=&-a_1^1e_1-a_1^2e_2;
%\end{array}
%$$
%and
%$$
%\begin{array}{rcl}
%T(e_2,e_3)&=&(\omega_3^1e_1+\omega_3^2e_2)(e_2)-(\omega_2^1e_1+\omega_2^3e_3)(e_3)-(a_2^1e_1+a_2^2e_2+a_2^3e_3)\\
%&=&\frac 1 2 e_1-\frac 1 2 e_1+a_2^3e_3-(a_2^1e_1+a_2^2e_2+a_2^3e_3)\\
%&=&-a_2^1e_1-a_2^2e_2.
%\end{array}
%$$
%\end{rmk}

\section{Surfaces in $M$}
Suppose that $S$ is an oriented differentiable two dimensional manifold in $M$. We get that $\dim(D\cap TS)\geq 1$, and as $de^3\wedge e^3=e^1\wedge e^2\wedge e^3$, the set where $\dim(D\cap TS)=2$ has empty interior. We denote by $\Sigma=\{x\in S:  \dim(D_x\cap T_xS)=2\}$ and by $S'=S-\Sigma$. The set $S'$ is open in $S$. In what follows we will suppose $\Sigma=\emptyset$, so  $S=S'$. With this hypothesis, the one dimensional vector subbundle $D\cap TS$  is well defined on $S$. Suppose $U\subset S$ is an open set such that we can define a unitary vector field $f_2$ with values in $D\cap TS$, so that $<f_2,f_2>=1$. Suppose $f_2=xe_1+ye_2$.
\begin{dfn}The unitary vector field $f_1\in \underline{D}$ defined by
$$f_1=ye_1-xe_2$$
is the \emph{horizontal normal} to $S$.
\end{dfn}
Then we can define $f^1\in \underline{(TM)^*|_S}$ by $f^1(f_1)=1$ and $f^1(TS)=0$. We call $f^1$ the \emph{horizontal conormal} to $S$. 
%\begin{dfn}The application
%$$\begin{array}{rcl}
%g=\exp\circ L^*\circ\eta^*:S&\rightarrow&\mathbb H^1\\
%p&\rightarrow&\exp(L_p^*(\eta^*(p)))
%\end{array}
%$$
%is the \emph{Gauss map} of $S$.
%\end{dfn}
If 
$$f_3=e_3-f^1(e_3)f_1,$$
then $\{f_2,f_3\}$ is a \emph{special} basis of $TS$ on the open set $U$. If $\alpha$ is the angle between $e_1$ and $f_1$, 
$$f_1=\cos\alpha\, e_1+\sin\alpha\, e_2,$$
where  $\alpha$ is a real function on $U$, reducing $U$ if necessary, then
$$f_2=-\sin\alpha\, e_1+\cos\alpha\, e_2,$$
and if we denote by $A=-f^1(e_3)$,
$$f_3=e_3+Af_1.$$
The dual basis of $(TM)^*$ on $S$ is 
\begin{equation*}%\label{dualf}
\begin{array}{lcl}
f^1&=&\cos\alpha\, e^1+\sin\alpha\, e^2-Ae^3,\\
f^2&=&-\sin\alpha\, e^1+\cos\alpha\, e^2,\\
f^3&=&e^3.
\end{array}
\end{equation*}
The inverse relations are
\begin{equation}\label{duale}
\begin{array}{lcl}
e^3&=&f^3,\\
e^1&=&\cos\alpha\,f^1-\sin\alpha\, f^2+A\cos\alpha\, f^3,\\
e^2&=&\sin\alpha\,f^1+\cos\alpha\, f^2+A\sin\alpha\, f^3.
\end{array}
\end{equation}

\section{The orthonormal basis of $TS$ in the scalar product $<,>_L$}
An orthonormal basis of $TS$ in $(M,<,>_L)$ is given by \begin{equation}\label{X2X3}
X_2^L=f_2 \mbox{ and }X_3^L=\frac{1}{\sqrt{L+A^2}}f_3=\sin\beta f_1+\cos\beta e_3^L,
\end{equation} 
where $\cos\beta=\frac{\sqrt{L}}{\sqrt{L+A^2}}$ and $\sin\beta=\frac{A}{\sqrt{L+A^2}}$. 
The normal vector in $<,>_L$ to $S$ is 
\begin{equation}\label{X1}
X_1^L=\frac{\sqrt{L}}{\sqrt{L+A^2}}f_1-\frac{A}{\sqrt{L+A^2}}e_3^L%=\frac{\sqrt{L}}{\sqrt{L+A^2}}f_1-\frac{A}{\sqrt{L+A^2}}(\frac{1}{\sqrt{L}}(f_3-Af_1))
=\cos\beta f_1-\sin\beta e_3^L.
\end{equation} 
%When $L\rightarrow +\infty$ it is clear thar $N^L\rightarrow f_1$, which is in  $S^1=\{(x,y,0):x^2+y^2=1\}$. This is not interesting. But the interesting notion is the dual one. Consider the dual 1-form 
%$(N^L)^*$ defined by $(N^L)^*(N^L)=1$ and $(N^L)^*(TS)=0$. A simple calculation shows that
%$$
%(N^L)^*=\frac{\sqrt{L}}{\sqrt{L+A^2}}f^1
%$$
%so $\lim_{L\rightarrow \infty}(N^L)^*=f^1$. We identify the image of $f^1$  in the cylinder $C=\{(x,y,z)\in\mathbb H^1:x^2+y^2=1\}$, defining a Gauss map by 
%$$
%g:p\in S\rightarrow (\cos\alpha,\sin\alpha, -A)\in C.
%$$
We can write the orthonormal basis as 
$$
\begin{array}{rcl}
X_1^L&=&\cos\beta \cos\alpha e_1+\cos\beta\sin\alpha e_2-\sin\beta e_3^L\\
X_2^L&=&-\sin\alpha e_1+\cos\alpha e_2\\
X_3^L&=&\sin\beta\cos\alpha e_1+\sin\beta\sin\alpha e_2+\cos\beta e_3^L
\end{array}
$$
and
$$
\begin{array}{rcl}
e^1&=&\cos\beta \cos\alpha X^1_L-\sin\alpha X^2_L+\sin\beta\cos\alpha X^3_L\\
e^2&=&\cos\beta\sin\alpha X^1_L+\cos\alpha X_L^2+\sin\beta\sin\alpha X^3_L\\
e^3_L&=&-\sin\beta X^1_L+\cos\beta X_L^3.
\end{array}
$$
As $d\sin\beta=\cos\beta d\beta$, we get
$$
d\beta=\frac{\sqrt{L}}{L+A^2}dA.
$$
Also
$$
\begin{array}{rcl}
\overline\nabla^Lf_1&=&d\alpha(-\sin\alpha e_1+\cos\alpha e_2)+\cos\alpha\Sigma_{j=1}^3\omega_1^{Lj}e_j^L+\sin\alpha\Sigma_{j=1}^3\omega_2^{Lj}e_j^L\\
&=&(d\alpha +\omega_1^{L2})f_2+(\cos\alpha\omega_1^{L3}+\sin\alpha\omega_2^{L3})e_3^L,
\\\end{array}
$$
\begin{equation}\label{X2}
\begin{array}{rcl}
\overline\nabla^Lf_2&=&d\alpha(-\cos\alpha e_1-\sin\alpha e_2)-\sin\alpha\Sigma_{j=1}^3\omega_1^{Lj}e_j^L+\cos\alpha\Sigma_{j=1}^3\omega_2^{Lj}e_j^L\\
&=&-(d\alpha +\omega_1^{L2})f_1+(-\sin\alpha\omega_1^{L3}+\cos\alpha\omega_2^{L3})e_3^L,
\\\end{array}
\end{equation}
and
$$
\begin{array}{rcl}
\overline\nabla^Le_3^L&=&\omega_3^{L1}e_1^L+\omega_3^{L2}e_2^L\\
&=&\omega_3^{L1}(\cos\alpha f_1-\sin\alpha f_2)+\omega_3^{L2}(\sin\alpha f_1+\cos\alpha f_2)\\
&=&(\cos\alpha \omega_3^{L1}+\sin\alpha \omega_3^{L2}) f_1+(-\sin\alpha \omega_3^{L1}+\cos\alpha \omega_3^{L2}) f_2.
\\\end{array}
$$
\section{The projection $\nabla^L$ of $\overline\nabla^L$ on $TS$}
The connection $\nabla^L$ on $TS$ is defined by
$$\nabla^L_XY=\overline\nabla^L_XY-<\overline\nabla^L_XY,X_1^L>_LX_1^L,$$
for $X,Y$ sections of $TS$.
We have 
$$\nabla^L X^L_2=<\overline\nabla^LX^L_2,X^L_3>_LX_3^L.$$
Taking into account \ref{X2} we get
\begin{equation}\label{Omega}
\Omega_2^{L3}:=<\overline\nabla^LX^L_2,X^L_3>_L=-\sin\beta(d\alpha+\omega_1^{L2})+\cos\beta(-\sin\alpha \omega_1^{L3}+\cos\alpha \omega_2^{L3}).
\end{equation}
As $X^1_L$ is null on $TS$ we get
$$
\nabla^L X_2^L=\Omega_2^{L3}%\otimes 
X_3^L.
$$
In the same way as $<\overline\nabla^LX^L_3,X^L_2>_L=-<\overline\nabla^LX^L_2,X^L_3>_L$, we obtain
$$\nabla^L X_3^L=-\Omega_2^{L3}
%\otimes 
X_2^L.$$

\section{The limit $K$ of curvatures $K^L$ of $S$}\label{curvgauss}

Now we will calculate the Gaussian curvature 
$$K^L=<R^L(X_2^L,X_3^L)X^L_3,X_2^L>_L,$$ 
where $R^L(X,Y)Z=\nabla_X\nabla_YZ-\nabla_Y\nabla_XZ-\nabla_{[X,Y]}Z$. Therefore
$$
\begin{array}{rcl}
<\nabla_{X_2^L}\nabla_{X_3^L}X_3^L,X_2^L>_L&=&<\nabla_{X_2^L}(\Omega^{L2}_3(X_3^L)X_2^L),X_2^L>_L\\
&=&<X_2^L(\Omega^{L2}_3(X_3^L))X_2^L+\Omega^{L2}_3(X_3^L)\nabla_{X_2^L}X_2^L),X_2^L>_L\\
&=&X_2^L(\Omega^{L2}_3(X_3^L)).
\end{array}
$$
In the same way
$$
<\nabla_{X^L_3}\nabla_{X_2^L}X_3^L,X_2^L>_L=X_3^L(\Omega^{L2}_3(X_2^L)), 
$$
and
$$
<\nabla_{[X^L_2,X^L_3]}X_3^L,X_2^L>_L=\Omega^{L2}_3([X^L_2,X^L_3])
$$
Then we have proved
\begin{prop}\label{KaL}
$K^L=d\Omega^{L2}_3(X^L_2,X^L_3).$
\end{prop}
\begin{dfn}
Define $K$ by
$$
K=\lim_{L\rightarrow \infty}K^L,
$$
if this limit exists.
\end{dfn}
It follows from proposition \ref{KaL} that
\begin{equation*}%\label{K}
K=%\lim_{L\rightarrow \infty}d\Omega^{L2}_3(X^L_2,X^L_3)=
\lim_{L\rightarrow \infty}d\Omega^{L2}_3(f_2,\frac{1}{\sqrt{L+A^2}}f_3)=\lim_{L\rightarrow \infty}\frac{1}{\sqrt{L}}d\Omega^{L2}_3(f_2,f_3)%=d\lim_{L\rightarrow \infty}\Omega^{L2}_3(f_2,\frac{f_3}{\sqrt{L}}).
\end{equation*}

%\Omega_2^{L3}:=<\overline\nabla^LX^L_2,X^L_3>&=&-\sin\beta(d\alpha+\omega_1^{L2})+\cos\beta(-\sin\alpha \omega_1^{L3}+\cos\alpha \omega_2^{L3})\\

\begin{prop}\label{dO}
%If the coefficients of forms $\omega_1^{L2}$, $\omega_1^{L3}$ e $\omega_2^{L3}$ in the basis $e^1$, $e^2$ and $e^3$ are rational functions in $\sqrt{L}$, and if $\lim_{L\rightarrow \infty} \frac{1}{\sqrt{L}}\omega_i^{L3}$   for $i=1,2$ and  $\lim_{L\rightarrow \infty} \frac{1}{L}\omega_1^{L2}$ exist, then 
We have
$$K=d\lim_{L\rightarrow \infty}\frac{1}{\sqrt{L}}\Omega^{L2}_3(f_2,f_3).$$
\end{prop}
\Pf We know from \ref{Omega} that
$$
\begin{array}{rcl}
\frac{1}{\sqrt{L}}\Omega_2^{L3}&=&-\frac{A}{\sqrt{L}\sqrt{L+A^2}}(d\alpha+\omega_1^{L2})+\frac{\sqrt{L}}{\sqrt{L}\sqrt{L+A^2}}(-\sin\alpha \omega_1^{L3}+\cos\alpha \omega_2^{L3})\\
&=&\frac{\sqrt{L}}{\sqrt{L+A^2}}\left(-\frac{A}{L}(d\alpha+\omega_1^{L2})+\frac{1}{\sqrt{L}}(-\sin\alpha \omega_1^{L3}+\cos\alpha \omega_2^{L3})\right)\\
&=&\frac{\sqrt{L}}{\sqrt{L+A^2}}\Theta_2^{L3},
\end{array}
$$
where 
$$\Theta_2^{L3}=-\frac{A}{L}(d\alpha+\omega_1^{L2})+\frac{1}{\sqrt{L}}(-\sin\alpha \omega_1^{L3}+\cos\alpha \omega_2^{L3}).$$
It follows that
$$
d\frac{1}{\sqrt{L}}\Omega_2^{L3}=d\frac{\sqrt{L}}{\sqrt{L+A^2}}\wedge \Theta_2^{L3}+\frac{\sqrt{L}}{\sqrt{L+A^2}}d\Theta_2^{L3}=\frac{-A\sqrt{L}}{\sqrt{L+A^2}^3}dA\wedge \Theta_2^{L3}+\frac{\sqrt{L}}{\sqrt{L+A^2}}d\Theta_2^{L3},
$$
and as we know by proposition \ref{limomega}  that 
$$
\begin{array}{rcl}
\lim_{L\rightarrow \infty} \Theta_2^{L3}&=&-\lim_{L\rightarrow \infty}\frac{A}{L}(d\alpha+\omega_1^{L2})+\lim_{L\rightarrow \infty}\frac{1}{\sqrt{L}}(-\sin\alpha \omega_1^{L3}+\cos\alpha \omega_2^{L3})\\
&=&-A\lim_{L\rightarrow \infty}\frac{1}{L}\omega_1^{L2}-\sin\alpha \lim_{L\rightarrow \infty}\frac{1}{\sqrt{L}}\omega_1^{L3}+\cos\alpha \lim_{L\rightarrow \infty}\frac{1}{\sqrt{L}}\omega_2^{L3}\\
&=&-A\omega_1^{2}-\sin\alpha \omega_1^{3}+\cos\alpha \omega_2^{L3}
\end{array}
$$
exists, we get
$$
\lim_{L\rightarrow \infty}d\frac{1}{\sqrt{L}}\Omega_2^{L3}=\lim_{L\rightarrow \infty}d\Theta_2^{L3},
$$
and by expressions \ref{omegaL} we get  that 
$$\lim_{L\rightarrow \infty}d\Theta_2^{L3}=d\lim_{L\rightarrow \infty}\Theta_2^{L3},$$
what ends the proof.\EPf
\begin{dfn}\label{principal}
We define 
$\Omega_2^3=\lim_{L\rightarrow \infty}\frac{1}{\sqrt{L}}\Omega_2^{L3}=-A\omega_1^{2}-\sin\alpha \omega_1^{3}+\cos\alpha \omega_2^{3}$. 
\end{dfn}
\begin{prop}
$\Omega_2^3=A f^3$.
\end{prop}
\Pf We get from definition \ref{limomega} and \ref{duale}
$$
\begin{array}{rcl}
\Omega_2^3&=&-A\omega_1^{2}-\sin\alpha \omega_1^{3}+\cos\alpha \omega_2^{3}\\
%&=&A\frac{1}{2}e^3+\sin\alpha(\frac{1}{2}e^2)+\cos\alpha(\frac{1}{2}e^1)\\
&=&\frac{1}{2}(\cos\alpha e^1+\sin\alpha e^2)+(\frac{1}{2}A) e^3\\
&=&A f^3.
\end{array}
$$
\EPf

\begin{prop}
$$
K=d\Omega_3^2(f_2,f_3)=-dA(f_2)-A^2.
$$
\end{prop}
\Pf It follows from \ref{duale} that
$$
\begin{array}{lcl}
df^3&=&de^3\\
&=&-e^1\wedge e^2\\
&=&-(-\sin\alpha\, f^2+A\cos\alpha\, f^3)\wedge (\cos\alpha\, f^2+A\sin\alpha\, f^3)\\
%&&-a_1^3(-\sin\alpha\, f^2)\wedge f^3-a_2^3(\cos\alpha\, f^2)\wedge f^3\\
&=&Af^2\wedge f^3$$
\end{array}
$$
then
$$
d\Omega_2^3=d(Af^3)=dA\wedge f^3+A^2f^2\wedge f^3,
$$
and the proposition follows. \EPf

\begin{rmk} In  \cite{BTV} and  \cite{WW} the calculus of $K$ was done applying   Gauss equation $K^L=\overline K^L+II^L$. Both terms $\overline K^L$ and $II^L$ are divergent, but the divergences cancel. In fact, if 
$$\overline\nabla^L X_1^L=\Omega_1^{L2}X_2^L+\Omega_1^{L3}X_3^L ,$$
we get 
$$
II^L(X_2^L,X_3^L)=\Omega_1^{L2}(X_2^L)\Omega_1^{L3}(X_3^L)-\Omega_1^{L2}(X_3^L)\Omega_1^{L3}(X_2^L)=\Omega_1^{L2}\wedge\Omega_1^{L3}(X_2^L,X_3^L),
$$
where 
$$\Omega_1^{L2}=-d\beta+\cos\beta(d\alpha+\omega_1^{L2})+(\sin\alpha\omega_1^{L3}-\cos\alpha\omega_2^{L3})$$
and
$$
\Omega_1^{L3}=(\cos\alpha\omega_1^{L3}+\sin\alpha\omega_2^{L3}).
$$
For another side 
$$
\overline K^L=<\overline R^L(X_2^L,X_3^L)X_3^L,X_2^L>_L=(d\Omega_3^{2L}-\Omega_1^{L2}\wedge\Omega_1^{L3})(X_2^L,X_3^L).
$$
It is straight verification  that $\lim_{L\rightarrow\infty}II^L$ diverge.
\end{rmk}

\section{The limit $k_n$ of normal curvatures $k_n^L$ of transverse curves in $S$}\label{curvnormal}

Suppose $\gamma(t)$ is a curve in $S$ such that $\gamma'(t)=x(t)f_2(\gamma(t))+y(t)f_3(\gamma(t))$, where $y(t)\neq 0$ for every $t$. Then $\gamma'(t)=x(t)X_2^L+y(t)\sqrt{L+A^2}X_3^L$ and the unitary tangent vector in the metric $<,>_L$ is given by the function 
$$
T^L=\frac{1}{\sqrt{x^2+y^2(L+A^2)}}(xX_2^L+y\sqrt{L+A^2}X_3^L).
$$
Let's write $x^L=\frac{x}{\sqrt{x^2+y^2(L+A^2)}}$ and $y^L=\frac{y\sqrt{L+A^2}}{\sqrt{x^2+y^2(L+A^2)}}$, so that 
$T^L=x^LX_2^L+y^LX_3^L$. Then
$$
\begin{array}{rcl}
\nabla^L_{T^L}T^L&=&\frac{d}{dt}x^LX_2^L+\frac{d}{dt}y^LX_3^L+x^L\nabla^L_{T^L}X_2^L+y^L\nabla^L_{T^L}X_3^L\vspace{.1cm}\\
&=&\frac{d}{dt}x^LX_2^L+\frac{d}{dt}y^LX_3^L+\Omega_2^{L3}(T^L)(-y^LX_2^L+x^LX_3^L)\\
\end{array}
$$
The normal vector to $T^L$ in $TS$ is  $N^L=-y^LX_2^L+x^LX_3^L$. Then $k_n^L=<\nabla^L_{T^L}T^L,N^L>$, so

\begin{equation}\label{kn}
k_n^L=-y_L\frac{d}{dt}x^L+x_L\frac{d}{dt}y^L+\Omega_2^{L3}(T^L)\vspace{.1cm}\\
%&=&-b_L\frac{d}{dt}a^L+a_L\frac{d}{dt}b^L+(-\sin\beta d\alpha(T^L)+\frac{\sqrt{L}}{2}\sin(2\beta)X^3_L(T^L))\\
%&=&\frac{1}{\sqrt{L+A^2}(a^2+b^2(L+A^2))}(abA\frac{d}{dt}A+(a\frac{d}{dt}b-b\frac{d}{dt}a)(L+A^2))\vspace{.1cm}\\
%&&-\sin\beta(a_Ld\alpha(X_2^L)+b_Ld\alpha(X_3^L))+\sqrt{L}\sin\beta\cos\beta b_L\vspace{.1cm}\\
%&=&\frac{1}{\sqrt{L+A^2}(a^2+b^2(L+A^2))}(abA\frac{d}{dt}A+(a\frac{d}{dt}b-b\frac{d}{dt}a)(L+A^2))\vspace{.1cm}\\
%&&-\frac{A}{\sqrt{L+A^2}}\frac{a}{\sqrt{a^2+b^2(L+A^2)}}d\alpha(f_2)-\frac{A}{\sqrt{L+A^2}}\frac{b\sqrt{L+A^2}}{\sqrt{a^2+b^2(L+A^2)}}d\alpha(\frac{1}{\sqrt{L+A^2}}f_3)\vspace{.1cm}\\
%&&+\sqrt{L}\frac{A}{\sqrt{L+A^2}}\frac{\sqrt{L}}{\sqrt{L+A^2}}\frac{b\sqrt{L+A^2}}{\sqrt{a^2+b^2(L+A^2)}}
\end{equation}
Observe that 
$$
\frac{d}{dt}x^L=%\frac{d}{dt}\frac{x}{\sqrt{x^2+y^2(L+A^2)}}=
\frac{dx}{dt}\frac{1}{\sqrt{x^2+y^2(L+A^2)}}-x\frac{x\frac{dx}{dt}+y\frac{dy}{dt}(L+A^2)+y^2A\frac{dA}{dt}}{\sqrt{x^2+y^2(L+A^2)}^3}
$$
so from similar formula for $\frac{d}{dt}y^L$, we get 
\begin{equation}\label{classical}
\lim_{L\rightarrow \infty}(-y_L\frac{d}{dt}x^L+x_L\frac{d}{dt}y^L)=0.
\end{equation}
Now
\begin{equation}\label{cerne}
\begin{array}{rcl}
\Omega_2^{L3}(T^L)&=&\frac{1}{\sqrt{x^2+y^2(L+A^2)}}\Omega_2^{L3}(xf_2+yf_3)%\\
%&=&\frac{1}{\sqrt{x^2+y^2(L+A^2)}}\frac{1}{\sqrt{L+A^2}}[-A(d\alpha+\omega_1^{L2})+\sqrt{L}(-\sin\alpha \omega_1^{L3}+\cos\alpha \omega_2^{L3})](xf_2+yf_3)
\end{array}
\end{equation}
\begin{dfn}
$$
k_n=\lim_{L\rightarrow\infty} k_n^L.
$$
\end{dfn}
\begin{prop}\label{knormal}
We have 
$$k_n=\frac{y}{|y|}A.$$
\end{prop}

\Pf It follows from definiton \ref{principal}, \ref{kn}, \ref{classical} and \ref{cerne} that
$$
\begin{array}{rcl}
k_n&=&\frac{1}{|y|}\Omega_2^3(xf_2+yf_3)\\
&=&\frac{1}{|y|}A f^3(xf_2+yf_3)\\
&=&\frac{y}{|y|}A
\end{array}
$$
\EPf

\section{The limit of Riemannian area element of $S$}\label{liRi}
It follows from \ref{X1} and \ref{X2X3} that 
%$$X_1^L=\frac{\sqrt{L}}{\sqrt{L+A^2}}f_1-\frac{A}{\sqrt{L+A^2}}e_3^L=\frac{\sqrt{L}}{\sqrt{L+A^2}}f_1-\frac{A}{\sqrt{L+A^2}}(\frac{1}{\sqrt{L}}(f_3-Af_1))=\frac{\sqrt{L+A^2}}{\sqrt{L}}f_1-\frac{A}{\sqrt{L+A^2}\sqrt{L}}f_3,$$
%$$X_2^L=f_2,$$ 
%$$X_3^L=\frac{1}{\sqrt{L+A^2}}f_3,$$ we obtain 
$$X^1_L=\frac{\sqrt{L}}{\sqrt{L+A^2}}f^1,\mbox{ } X^2_L=f^2 \mbox{ and }X^3_L=\sqrt{L+A^2}f^3+\frac{A}{\sqrt{L+A^2}}f^1.$$ 
Therefore on $S$ we get
$$
d\sigma_L=X^2_L\wedge X^3_L=\sqrt{L+A^2}f^2\wedge f^3%+\frac{A}{\sqrt{L+A^2}}f^2\wedge f^1.
$$
since that $f^1$ is null on $TS$. We can see that $\lim_{L\rightarrow \infty}K^Ld\sigma_L$ does not exist. In \cite{BTV} and \cite{WW}, to get an area form on $S$ it was necessary to multiply $d\sigma_L$ by $\frac{1}{\sqrt{L}}$ and take the limit as $L$ goes to infinity to obtain a surface form, i.e.,
$$
d\sigma=\lim_{L\rightarrow \infty }\frac{1}{\sqrt{L}}d\sigma_L=f^2\wedge f^3,
$$
which is  the Hausdorff measure on $S$. Therefore
%$$
%K^Ld\sigma_L=\left(\frac{L}{(L+A^2)^2} d\alpha\wedge dA(f_3,f_2)-\frac{L^2}{(L+A^2)^2}dA(f_2)-\frac{L}{L+A^2}A^2\right)\left(\sqrt{L+A^2}f^2\wedge f^3+\frac{A}{\sqrt{L+A^2}}f^2\wedge f^1 \right)
%$$
%sot
\begin{equation}\label{KL}
\lim_{L\rightarrow \infty }\frac{1}{\sqrt{L}}K^Ld\sigma_L=Kf^2\wedge f^3=Kd\sigma.
\end{equation}

%$$
%K^Ld\sigma_L=\frac{L}{\sqrt{L+A^2}} \left(\frac{1}{L+A^2} (d\alpha\wedge dA(f_3,f_2)-LdA(f_2))-A^2\right)\left(f^2\wedge f^3+\frac{A}{L+A^2}f^2\wedge f^1 \right).
%$$

\section{The limit of length elements}

%The normal in $TS$ to the tangent vector $T^L$ is $N^L=-b^LX_2^L+a^LX_3^L$, therefore 
The length element in the metric $<,>_L$ on $\gamma$ is
$$
ds_L=%-i(N^L)d\sigma_L=-i(N^L)X_L^2\wedge X_L^3=
x^LX^2_L+y^LX^3_L.
$$
%But $\lim_{L\rightarrow\infty}\sin\beta=0$, $\lim_{L\rightarrow\infty}\cos\beta=1$, so  it follows from $X^2_L=\sec\alpha e^2-\tan\alpha\tan\beta e^3_L$
%and $X^3_L=\sec\beta e^3_L$ that $\lim _{L\rightarrow\infty}X^2_L=\sec\alpha e^2 $ and $\lim _{L\rightarrow\infty}X^3_L=\lim _{L\rightarrow\infty}\sqrt{L+A^2}e^3$ does not exist. 
As $\lim _{L\rightarrow\infty}x^L=0$ and $\lim _{L\rightarrow\infty}y^L=\frac{y}{|y|}$, we get
$$
\lim _{L\rightarrow\infty}ds_L=\lim _{L\rightarrow\infty}(x^Lf^2+y^L\sqrt{L+A^2}f^3)=\frac{y}{|y|}f^3\lim _{L\rightarrow\infty}\sqrt{L+A^2}
$$
that does not exist. But as  in \cite{BTV} and section \ref{liRi}, if we multiply by $\frac{1}{\sqrt{L}}$ we obtain
$$
ds=\lim _{L\rightarrow\infty}\frac{1}{\sqrt{L}}ds_L=\frac{y}{|y|}f^3,
$$
which is Hausdorff measure for transversal curves. It follows that
\begin{equation}\label{kL}
\lim _{L\rightarrow\infty}\frac{1}{\sqrt{L}}k_n^Lds_L=k_nds=Af^3=\Omega_2^3.
\end{equation}

\section{The Gauss-Bonnet theorem}
The proof of Gauss-Bonnet theorem in \cite{BTV} and \cite{WW} was done taking limits of  Gauss-Bonnet formulas in $(R^3,<,>_L)$ 
$$
\int_S\frac{1}{ \sqrt{L}}K^L d\sigma_L+\int_{\partial S}\frac{1}{ \sqrt{L}}k_n^L ds_L=\frac{1}{ \sqrt{L}}2\pi\chi(S)
$$
as $L$ goes to infinity to obtain $\int_SK^\infty d\sigma+\int_{\partial S}k_n ds=0$.

We will give below a straightforward proof due to  expressions of $K$ and $k_n$ obtained in sections \ref{curvgauss} and \ref{curvnormal}. We will restrict our theorem to regions where points are non singular and the boundary is constituted by transverse curves.

%\begin{thm}
%Let $S$ be a differentiable compact surface in $\mathbb H^1$ without characteristic points, and with finitely many boundary components $\partial S_i$, $i=1,\cdots,n$, given by smooth closed transverse curves $\gamma_i$. Then
%$$\int_S K^\infty d\sigma+\Sigma_{i=1}^n\int_{\gamma_i}k_nds=0.$$
%\end{thm}
%\Pf We will omit the technical details of the proof, since that they are well known:
%$$\int_S K^\infty d\sigma=\int_S (-dA(f_2)-A^2)f^2\wedge f^3=$$
%K^\infty d\sigma
%The theorem  follows from (\ref{KL}) and (\ref{kL}) taking limit as $L$ tends to infinity of Gauss-Bonnet theorem in $(\mathbb R^3,g_L)$:
%$$
%\int_S\frac{1}{ \sqrt{L}}K^L d\sigma_L+\Sigma_{i=1}^n\int_{\gamma_i}\frac{1}{ \sqrt{L}}k_n^L ds_L=\frac{1}{ \sqrt{L}}2\pi\chi(S).
%$$
%%%%%%%%
%%%%%%%%%
Let be $R\subset S$ a fundamental set, and $c$ a fundamental $2$-chain such that $|c|=R$. The oriented curve $\gamma=\partial c$ is the bounding curve of $R$. The curve $\gamma$ is piecewise differentiable, and composed by differentiable curves $\gamma_j:[s_j,s_{j+1}]\rightarrow S$, $j=1,\ldots,r$, with $\gamma_1(s_1)=\gamma_r(s_{r+1})$  and $\gamma_j(s_{j+1})=\gamma_{j+1}(s_{j+1})$, for $j=1,\ldots,r-1$. 
\begin{thm}(Gauss-Bonnet formula) Let ${R}$ be contained in a coordinate domain $U$ of $S$, let the bounding curve $\gamma$ of $R$ be a simple closed transverse curve. Then
$$\int_{\gamma}k_n+\int_R K=0,$$
where $k_n=\lim_{L\rightarrow\infty} k_n^L$  on $\gamma$ and $K=\lim_{L\rightarrow \infty}K^L$  on $R$.
\end{thm}
\Pf %Let $\gamma_1,\ldots,\gamma_r$ be the $C^{\infty}$ pieces of $\gamma$ with $\gamma_j$ defined on the interval $[s_j,s_{j+1}]$, and $\gamma_j(s_{j+1})=\gamma_{j+1}(s_{j+1})$, for $j=1,\ldots,r-1$, and $\gamma_r(s_{r+1})=\gamma_1(s_1)$.
 From Stokes theorem and using \ref{KL} and \ref{kL} we get
%In each $C^\infty$ piece of $\gamma$ we have the positive orientation $T$ and the curvature $\nabla_TT=\epsilon kf_1$. Then from (\ref{K}), Propositions \ref{pnabla} \ref{fork} and
%$$f_2=\frac{1}{f^2(\gamma_j')}\gamma_j'-\frac{f^1(\gamma_j')}{f^2(\gamma_j')}f_1=-\epsilon T-\frac{f^1(\gamma_j')}{f^2(\gamma_j')}f_1,$$
%since that $\epsilon=-\frac{|f^2(\gamma_j')|}{f^2(\gamma_j')}$, we obtain
$$
\int_RK=\int_cK d\sigma=\int_cK f^2\wedge f^3=%\int_c d\Omega^3_2(f_2,f_3) f^2\wedge f^3=
\int_c d\Omega_3^2=-\int_{\partial c}\Omega_2^3=-\int_{\partial c}k_nds=-\int_{\gamma}k_n.
$$
\EPf
%%%%%%%%
%%%%%%%%%

\end{document}